\documentclass[letterpaper]{amsart}
\pdfoutput=1
\newtheorem{theorem}{Theorem}[section]

\newtheorem{proposition}[theorem]{Proposition}
\newtheorem{corollary}[theorem]{Corollary}

\theoremstyle{definition}

\theoremstyle{remark}

\usepackage{hyperref}
\usepackage[top=1.0in, left=1.0in, right=1.0in]{geometry}
\usepackage[all,cmtip]{xy}
\usepackage{xypic}
\usepackage{tikz,tikz-cd}
\usepackage{pgf}
\usetikzlibrary{positioning}
\usetikzlibrary{calc}
\usetikzlibrary{arrows}
\usetikzlibrary{shapes.multipart}
\usepackage{graphicx}
\usepackage{xcolor}
\definecolor{darkgreen}{rgb}{0.33, 0.55, 0.13}
\definecolor{darkpink}{rgb}{0.91, 0.33, 0.5}
\definecolor{electriccyan}{rgb}{0.0, 1.0, 1.0}
\definecolor{electricultramarine}{rgb}{0.25, 0.0, 1.0}
\definecolor{greenyellow}{rgb}{0.30, 0.82, 0.0}
\definecolor{lightbrown}{rgb}{0.71, 0.4, 0.11}
\definecolor{prune}{rgb}{0.45, 0.11, 0.11}
\usepackage{color}
\usepackage{amsmath}
\usepackage{amsfonts}
\usepackage{euler}
\usepackage{amssymb}
\usepackage{epsfig}
\usepackage{rotating}
\usepackage{graphics}
\newcommand{\be}{\begin{equation}}

\newcommand{\ee}{\end{equation}}

\newcommand{\dz}{\wedge}

\newcommand{\ba}{\begin{array}}

\newcommand{\ea}{\end{array}}

\newcommand{\beq}{\begin{eqnarray}}

\newcommand{\eeq}{\end{eqnarray}}

\newtheorem{lm}{lemma}

\newtheorem{thee}{theorem}

\newtheorem{proo}{proposition}

\newtheorem{co}{corollary}

\newtheorem{rem}{remark}

\newtheorem{deff}{definition}

\newcommand{\bd}{\begin{deff}}

\newcommand{\ed}{\end{deff}}

\newcommand{\bl}{\begin{lm}}

\newcommand{\el}{\end{lm}}

\newcommand{\bp}{\begin{proo}}

\newcommand{\ep}{\end{proo}}

\newcommand{\bt}{\begin{thee}}

\newcommand{\et}{\end{thee}}

\newcommand{\bc}{\begin{co}}

\newcommand{\ec}{\end{co}}

\newcommand{\brm}{\begin{rem}}

\newcommand{\erm}{\end{rem}}

\newcommand{\der}{{\rm d}}

\usepackage{t1enc}
\def\frak{\mathfrak}

\newcommand{\newc}{\newcommand}

\let\ccdot\cdot
\def\cdot{\hbox to 2.5pt{\hss$\ccdot$\hss}}

\newc{\aR}{\mbox{\boldmath{$ R$}}}
\newc{\aS}{\mbox{\boldmath{$ S$}}}
\newc{\aT}{\mbox{\boldmath{$ T$}}}
\newc{\aW}{\mbox{\boldmath{$ W$}}}

\newc{\aK}{\mbox{\boldmath{$ K$}}}
\newc{\aL}{\mbox{\boldmath{$ L$}}}

\usepackage{amssymb}
\usepackage{amscd}

\newcommand{\hook}{\raisebox{-0.35ex}{\makebox[0.6em][r]
{\scriptsize $-$}}\hspace{-0.15em}\raisebox{0.25ex}{\makebox[0.4em][l]{\tiny
 $|$}}}

\newcommand{\bma}{\begin{pmatrix}}
\newcommand{\ema}{\end{pmatrix}}

\newc{\obstrn}[2]{B^{#1}_{#2}}

\newcommand{\rpl}                         
{\mbox{$
\begin{picture}(12.7,8)(-.5,-1)
\put(0,0.2){$+$}
\put(4.2,2.8){\oval(8,8)[r]}
\end{picture}$}}

\newcommand{\lpl}                         
{\mbox{$
\begin{picture}(12.7,8)(-.5,-1)
\put(2,0.2){$+$}
\put(6.2,2.8){\oval(8,8)[l]}
\end{picture}$}}

\usepackage{ifthen}

\newc{\tensor}[1]{#1}
\newc{\Mvariable}[1]{\mbox{#1}}
\newc{\down}[1]{{}_{#1}}
\newc{\up}[1]{{}^{#1}}

\newc{\JulyStrut}{\rule{0mm}{6mm}}
\newc{\midtenPan}{\mbox{\sf S}}
\newc{\midten}{\mbox{\sf T}}
\newc{\midtenEi}{\mbox{\sf U}}
\newc{\ATen}{\mbox{\sf E}}
\newc{\BTen}{\mbox{\sf F}}
\newc{\CTen}{\mbox{\sf G}}

\def\sideremark#1{\ifvmode\leavevmode\fi\vadjust{\vbox to0pt{\vss
 \hbox to 0pt{\hskip\hsize\hskip1em
 \vbox{\hsize3cm\tiny\raggedright\pretolerance10000
 \noindent #1\hfill}\hss}\vbox to8pt{\vfil}\vss}}}%

\newcommand{\Span}{\mathrm{Span}}

\numberwithin{equation}{section}

\newcounter{romenumi}
\newcommand{\labelromenumi}{(\roman{romenumi})}

\begin{document}
\newcommand{\bbS}{\mathbb{S}}
\newcommand{\bbR}{\mathbb{R}}
\newcommand{\bbK}{\mathbb{K}}
\newcommand{\sog}{\mathbf{SO}}
\newcommand{\spg}{\mathbf{Sp}}
\newcommand{\glg}{\mathbf{GL}}
\newcommand{\slg}{\mathbf{SL}}
\newcommand{\og}{\mathbf{O}}
\newcommand{\soa}{\frak{so}}
\newcommand{\spa}{\frak{sp}}
\newcommand{\gla}{\frak{gl}}
\newcommand{\sla}{\frak{sl}}
\newcommand{\sua}{\frak{su}}
\newcommand{\sug}{\mathbf{SU}}
\newcommand{\cspg}{\mathbf{CSp}}
\newcommand{\gat}{\tilde{\gamma}}
\newcommand{\Gat}{\tilde{\Gamma}}
\newcommand{\thet}{\tilde{\theta}}
\newcommand{\Thet}{\tilde{T}}
\newcommand{\rt}{\tilde{r}}
\newcommand{\st}{\sqrt{3}}
\newcommand{\kat}{\tilde{\kappa}}
\newcommand{\kz}{{K^{{~}^{\hskip-3.1mm\circ}}}}
\newcommand{\bv}{{\bf v}}
\newcommand{\di}{{\rm div}}
\newcommand{\curl}{{\rm curl}}
\newcommand{\cs}{(M,{\rm T}^{1,0})}
\newcommand{\tn}{{\mathcal N}}
\newcommand{\ten}{{\Upsilon}}
\title{On certain classes of $\spg(2,\bbR)$ symmetric $G_2$ structures}
\vskip 1.truecm
\author{Pawe\l~ Nurowski} \address{Centrum Fizyki Teoretycznej,
Polska Akademia Nauk, Al. Lotnik\'ow 32/46, 02-668 Warszawa, Poland}
\email{nurowski@cft.edu.pl}
\thanks{Support: This work was supported by the Polish National Science Centre (NCN) via the grant number 2018/29/B/ST1/02583.}

\date{\today}
\begin{abstract}
We find two different families of $\spg(2,\bbR)$ symmetric $G_2$ structures in seven dimensions. These are $G_2$ structures with $G_2$ being the split real form of the simple exceptional complex Lie group $G_2$. The first family has $\tau_2\equiv 0$, while the second family has $\tau_1\equiv\tau_2\equiv 0$. The families are different in the sense that the first one lives on a homogoneous space $\spg(2,\bbR)/\slg(2,\bbR)_l$, and the second one lives on a homogeneous space $\spg(2,\bbR)/\slg(2,\bbR)_s$. Here $\slg(2,\bbR)_l$ is an $\slg(2,\bbR)$  corresponding to the $\sla(2,\bbR)$ related to the long roots in the root diagram of $\spa(2,\bbR)$, and  $\slg(2,\bbR)_s$ is an $\slg(2,\bbR)$  corresponding to the $\sla(2,\bbR)$ related to the short roots in the root diagram of $\spa(2,\bbR)$. 
\end{abstract}
\maketitle
\section{Introduction: a question of Maciej Dunajski}
Recently, together with C. D. Hill \cite{car}, we uncovered an $\spg(2,\bbR)$ symmetry of the nonholonomic kinematics of a car. I talked about this at the Abel Symposium in \r{A}lesund, Norway, this June. After my talk Maciej Dunajski, intrigued by the root diagram of $\spa(2,\bbR)$ which appeared at the talk, asked me if using it I can see a $G_2$ structure on a 7-dimensional homogeneous space $M=\spg(2,\bbR)/\slg(2,\bbR)$.\\
\centerline{\includegraphics[height=6cm]{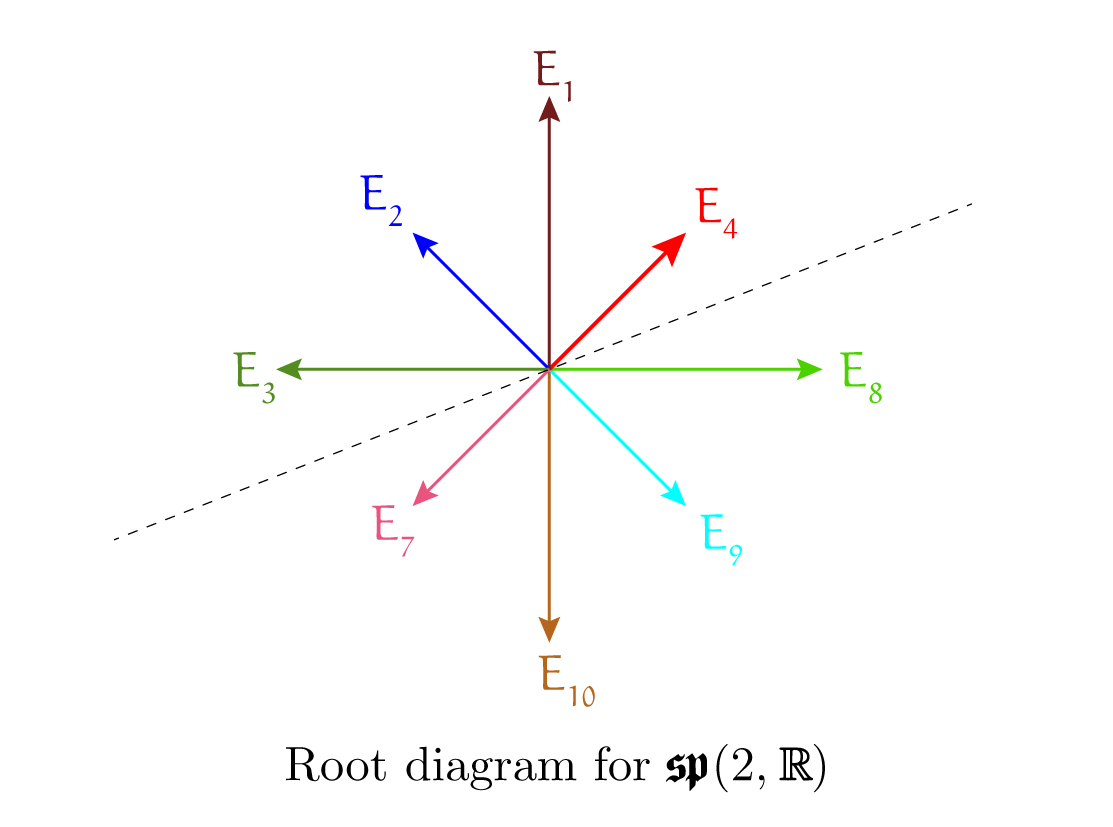}}
My immediate answer was: `I can think about it, but I have to know which of the $\slg(2,\bbR)$ subgroups of $\spg(2,\bbR)$ I shall use to built $M$'. The reason for the `but' word in my answer was that there are at least two $\slg(2,\bbR)$ subgroups of $\spg(2,\bbR)$, which lie quite differently in there. One can see them in the root diagram above: the first $\slg(2,\bbR)$ corresponds to the \emph{long} roots, as e.g. $E_1$ and $E_{10}$, whereas the second one corresponds to the \emph{short} roots, as e.g. $E_2$ and $E_9$. Since Maciej never told me which $\slg(2,\bbR)$ he wants, I decided to consider both of them, and to determine what kind of $G_2$ structures one can associate with each choice of subgroups, separately.

I emphasize that in the below considerations I will use the \emph{split real form} of the simple exceptional Lie group $G_2$. Therefore the corresponding $G_2$ structure metrics will \emph{not} be riemannian. They will have signature $(3,4)$.

\section{The Lie algebra $\spa(2,\bbR)$}
    The Lie algebra $\spa(2,\bbR)$ is given by the $4\times 4$ real matrices
    $$
    E=(E^\alpha{}_\beta)=
    \bma
    {a_5}&{a_7}&{a_9}&{2a_{10}}\\
    {-a_4}&{a_6}&{a_8}&{a_9}\\
    {a_2}&{a_3}&{-a_6}&{-a_7}\\
    {-2a_1}&{a_2}&{a_4}&{-a_5}
      \ema,
      $$
      where the coefficients $a_I$, $I=1,2,\dots 10$, are real constants. The commutator in $\spa(2,\bbR)$ is the usual commutator $[E,E']=E\cdot E'-E'\cdot E$ of two matrices $E$ and $E'$. We start with the following basis $(E_I)$, 
      $$E_I=\frac{\partial E}{\partial a_I},\quad I=1,2,\dots 10,$$
      in $\spa(2,\bbR)$.

     In this basis, modulo the antisymmetry, we have the following nonvanishing commutators: $[{E_1},{\color{black}E_5}]=2{E_1}$, $[{E_1},{E_7}]={-2E_2}$,  $[{E_1},{E_9}]={-2E_4}$, $[{E_1},{E_{10}}]={\color{black}4E_5}$, $[{E_2},{E_4}]={E_1}$, $[{E_2},{\color{black}E_5}]={E_2}$, $[{E_2},{\color{black}E_6}]={E_2}$, $[{E_2},{E_7}]={2E_3}$, $[{E_2},{E_8}]={E_4}$, $[{E_2},{E_9}]={\color{black}-E_5-E_6}$, $[{E_2},{E_{10}}]={-2E_7}$, $[{E_3},{E_4}]={-E_2}$, $[{E_3},{\color{black}E_6}]={2E_3}$, $[{E_3},{E_8}]={\color{black}-E_6}$, $[{E_3},{E_9}]={-E_7}$, $[{E_4},{\color{black}E_5}]={E_4}$, $[{E_4},{\color{black}E_6}]={-E_4}$, $[{E_4},{E_7}]={\color{black}E_5-E_6}$, $[{E_4},{E_9}]={-2E_8}$, $[{E_4},{E_{10}}]={-2E_9}$,  $[{\color{black}E_5},{E_7}]={E_7}$, $[{\color{black}E_5},{E_9}]={E_9}$, $[{\color{black}E_5},{E_{10}}]={2E_{10}}$, $[{\color{black}E_6},{E_7}]={-E_7}$, $[{\color{black}E_6},{E_8}]={2E_8}$, $[{\color{black}E_6},{E_9}]={E_9}$, $[{E_7},{E_8}]={E_9}$,  $[{E_7},{E_9}]={E_{10}}$.

     We see that there are at least \emph{two} $\sla(2,\bbR)$ Lie algebras here.
     The first one is
     $$\sla(2,\bbR)_l=\Span_\bbR(E_1,E_5,E_{10}),$$
     and the second is
     $$\sla(2,\bbR)_s=\Span_\bbR(E_2,E_5+E_6,E_9).$$
     The reason for distinguishing these two is as follows:

The eight 1-dimensional vector subspaces $\mathfrak{g}_I=\Span(E_I)$, $I=1,2,3,4,7,8,9,10$, of $\spa(2,\bbR)$ are the \emph{root spaces} of this Lie algebra. They correspond to the Cartan subalgebra of $\spa(2,\bbR)$ given by $\mathfrak{h}=\Span(E_5,E_6)$. It follows that the pairs $(E_I,E_J)$ of the root vectors, such that $I+J=11$, $I,J\neq 5,6$, correspond to the \emph{opposite roots} of $\sla(2,\bbR)$. Knowing the Killing form for $\sla(2,\bbR)$, which in the basis $(E_I)$, and its dual basis $(E^I)$, $E_I\hook E^J=\delta^J{}_I$, is         
$$K=\tfrac{1}{12}K_{IJ}E^I\odot E^J=-4{E^1} \odot {E^{10}}+2 {E^2}\odot {E^9}+{E^3}\odot {E^8} -2 {E^4}\odot {E^7}+{\color{black}E^5}\odot {\color{black}E^5}+{\color{black}E^6}\odot {\color{black}E^6},$$
one can see that the roots corresponding to the root vectors $(E_1,E_{10})$ and $(E_3,E_8)$ are \emph{long}, and the roots cooresponding to the root vectors $(E_2,E_9)$ and $(E_4,E_7)$ are \emph{short}. Thus the Lie algebra $\sla(2,\bbR)_l$ corresponding to the root vectors $(E_1,E_{10})$, and in turn to the \emph{long} roots, lies quite different in $\spa(2,\bbR)$ than the Lie algebra $\sla(2,\bbR)_s$ corresponding to the short roots associated with the root vectors $(E_2,E_9)$.  

\section{$G_2$ structures on $\spg(2,\bbR)/\slg(2,\bbR)_l$}
\subsection{Compatible pairs $(g,\phi)$ on $M_l$}
To consider the homogeneous space $M_l=\spg(2,\bbR)/\slg(2,\bbR)_l$ it is convenient to change the basis $(E_I)$ in $\spg(2,\bbR)$ to a new one, $(e_I)$, in which the last three vectors span $\slg(2,\bbR)_l$. Thus we take:
$$e_1=E_2,\,\,e_2=E_3,\,\,e_3=E_4,\,\,e_4=E_6,\,\,e_5=E_7,\,\,e_6=E_8,\,\,e_7=E_9,\,\,e_8=E_1,\,\,e_9=E_5,\,\,e_{10}=E_{10}.$$
If now, one considers $(e_I)$ as the basis of the Lie algebra of invariant vector fields on the Lie group $\spg(2,\bbR)$ then the dual basis $(e^I)$, $e_I\hook e^J=\delta^J{}_I$, of the left invariant forms on $\spg(2,\bbR)$ satisfies:
\be\begin{aligned}
  \der e^1=&-e^1\dz(e^4+e^9)+e^2\dz e^3-2e^5\dz e^8\\
    \der e^2=&-2e^1\dz e^5-2e^2\dz e^4\\
      \der e^3=&-e^1\dz e^6+e^3\dz(e^4-e^9)-2e^7\dz e^8\\
        \der e^4=&e^1\dz e^7+e^2\dz e^6+e^3\dz e^5\\
          \der e^5=&2e^1\dz e^{10}+e^2\dz e^7+e^5\dz (e^9-e^4)\\
            \der e^6=&2e^3\dz e^7-2e^4\dz e^6\\
              \der e^7=&2e^3\dz e^{10}-e^5\dz e^6+e^7\dz (e^4+e^9)\\
                \der e^8=&-e^1\dz e^3-2e^8\dz e^9\\
                  \der e^9=&e^1\dz e^7-e^3\dz e^5-4e^8\dz e^{10}\\
                    \der e^{10}=&-e^5\dz e^7-2e^9\dz e^{10}.
\end{aligned}\label{mc}\ee
Here we used the usual formula relating the structure constants $c^I{}_{JK}$, from $[e_J,e_K]=c^I{}_{JK}e_I$, to the differentials of the Maurer-Cartan forms $(e^I)$, $\der e^I=-\tfrac12 c^I{}_{JK}e^J\dz e^K$.

In this basis the Killing form on $\spg(2,\bbR)$ is
$$K=\tfrac{1}{12}c^I{}_{JK}c^K{}_{LI}e^J\odot e^L=(e^4)^2-2e^3\odot e^5+e^2\odot e^6+2e^1\odot e^7+(e^9)^2-4e^8\odot e^{10}.$$
Here we have used the notation $e^I\odot e^J=\tfrac12(e^I\otimes e^J+e^J\otimes e^I)$, $(e^I)^2=e^I\odot e^I$. 

Looking at the equations \eqref{mc} one sees that $\spg(2,\bbR)$ has the structure of the principal $\slg(2,\bbR)$ fiber bundle $\slg(2,\bbR)_l\to \spg(2,\bbR)\to M_l=\spg(2,\bbR)/\slg(2,\bbR)_l$ over the homogeneous space $M_l=\spg(2,\bbR)/\slg(2,\bbR)_l$.

Indeed, the 3-dimensional distribution $D_l$, generated by the vector fields $X$ on $\spg(2,\bbR)$ annihilating the span of the 1-forms $(e^1,e^2,\dots,e^7)$, is integrable, $\der e^\mu\dz e^1\dz e^2\dots \dz e^7\equiv 0$, $\mu=1,2\dots,7$,  so that we have a well defined 7-dimensional leaf space $M_l$ of the corresponding foliation. Moreover, the Maurer-Cartan equations \eqref{mc}, restricted to a leaf defined by $(e^1,e^2,\dots,e^7)\equiv 0$, reduce to $ \der e^8=-2e^8\dz e^9$,
$\der e^9=-4e^8\dz e^{10}$, $\der e^{10}=-2e^9\dz e^{10}$, showing that each leaf can be identified with the Lie group $\slg(2,\bbR)_l$. Thus the projection $\spg(2,\bbR)\to M_l$ from the Lie group $\spg(2,\bbR)$  to the leaf space $M_l$ is the projection to the homogeneous space $M_l=\spg(2,\bbR)/\slg(2,\bbR)_l$.  

From now on, in this Section,  I will use Greek indices $\mu,\nu$, etc., to run from 1 to 7. They number the first seven basis elements in the bases $(e_I)$ and $(e^I)$.

Now, I look for all bilinear symmetric forms $g=g_{\mu\nu}e^\mu\odot e^\nu$ on $\spg(2,\bbR)$, with constant coefficients $g_{\mu\nu}=g_{\nu\mu}$, which are constant along the leaves of the foliation defined by $D_l$. Technically, I search for those $g$ whose Lie derivative with respect to any vector field $X$ from $D_l$ vanishes, \be{\mathcal L}_Xg=0\,\, \mathrm{for\,\, all}\,\, X\,\, \mathrm{in}\,\, D_l.\label{c1}\ee
I have the following proposition:  
\begin{proposition}\label{pr1}
  The most general $g=g_{\mu\nu}e^\mu\odot e^\nu$ satisfying condition \eqref{c1} is
  $$g=g_{22}(e^2)^2+2g_{24}e^2\odot e^4+g_{44}(e^4)^2+2g_{35}(e^3\odot e^5-e^1\odot e^7)+2g_{26}e^2\odot e^6+2g_{46}e^4 \odot e^6+g_{66}(e^6)^2.$$
\end{proposition}
Thus I have a 7-parameter family of bilinear forms on $\spg(2,\bbR)$ that \emph{descend} to well defined pseudoriemannian metrics on the leaf space $M_l$. Note that, the restriction of the Killing form $K$ to the space where $(e^8,e^9,e^{10})\equiv 0$ is in this family. This corresponds to $g_{22}=g_{24}=g_{46}=0$ and $g_{44}=2g_{26}=-g_{35}=1$.

Since the aim of my note is \emph{not} to be exhaustive, but rather to show how to produce $G_2$ structures on $\spg(2,\bbR)$ homogeneous spaces, from now on I will restirict myself to only one $\slg(2,\bbR)_l$ invariant bilinear form $g$ on $\spg(2,\bbR)$, namely to
\be g_K=(e^4)^2-2e^3\odot e^5+e^2\odot e^6+2e^1\odot e^7,\label{metl}\ee
coming from the restriction of the Killing form. It follows from the Proposition \ref{pr1} that this form is a well defined $(3,4)$ signature metric on the quotient space $M_l=\spg(2,\bbR)/\slg(2,\bbR)_l$.

I now look for the 3-forms $\phi=\tfrac16\phi_{\mu\nu\rho}e^\mu\dz e^\nu\dz e^\rho$ on $\spg(2,\bbR)$ that are constant along the leaves of the distribution $D_l$, i.e. such that
\be{\mathcal L}_X\phi=0\,\, \mathrm{for\,\, all}\,\, X\,\, \mathrm{in}\,\, D_l.\label{c2}\ee
Then, I have the following proposition.
\begin{proposition}\label{pr2}
  There is a 10-parameter family of 3-forms $\phi=\tfrac16\phi_{\mu\nu\rho}e^\mu\dz e^\nu\dz e^\rho$ on $\spg(2,\bbR)$ which satisfy condition \eqref{c2}.
  The general formula for them is:
  $$\phi=f e^{125}+a(e^{235}-e^{127})+pe^{145}+q(e^{147}+e^{345})+se^{156}+t(e^{356}-e^{167})+he^{237}+be^{246}+re^{347}+ue^{367}.$$
  Here $e^{\mu\nu\rho}=e^\mu\dz e^\nu\dz e^\rho$, and $a$, $b$, $f$, $h$, $p$, $q$, $r$, $s$, $t$ and $u$ are real constants.
    \end{proposition}
Thus there is a 10-parameter family of 3-forms $\phi$ that descends from $\spg(2,\bbR)$ to the $\spg(2,\bbR)$ homogeneous space $M_l=\spg(2,\bbR)/\slg(2,\bbR)_l$, and is well defined there.

Now, I introduce an important notion of \emph{compatibility} of a pair $(g,\phi)$ where $g$ is a metric, and $\phi$ is a 3-form on a 7-dimensional oriented manifold $M$. The pair $(g,\phi)$ on $M$ \emph{is compatible} if and only if
$$(X\hook\phi)\dz(X\hook\phi)\dz\phi\,=\,3\,g(X,Y)\,\mathrm{vol}(g),\quad\quad \forall X,Y\in \mathrm{T}M.$$
Here $\mathrm{vol}(g)$ is a \emph{volume form} on $M$ related to the metric $g$. 

Restricting, as I did, to the $\spg(2,\bbR)$ invariant metric $g_K$ on $M_l$ as in \eqref{metl}, I now ask which of the 3-forms $\phi$ from Proposition \ref{pr2} are compatible with the metric \eqref{metl}. In other words, I now look for the constants  $a$, $b$, $f$, $h$, $p$, $q$, $r$, $s$, $t$ and $u$ such that
\be(e_\mu\hook\phi)\dz(e_\nu\hook\phi)\dz\phi\,=\,3\,g_K(e_\mu,e_\nu)\,e^1\dz e^2\dz e^3\dz e^4\dz e^5\dz e^6\dz e^7,\label{coml}\ee
for $g=g_K$ given in \eqref{metl}.

I have the following proposition.
\begin{proposition}
  The general solution to the equations \eqref{coml} is given by
  $$b=\tfrac12,\,\,f=\frac{ap}{1-q},\,\,h=\frac{a(q-1)}{p},\,\,r=\frac{q^2-1}{p},\,\,s=\frac{p(1-q)}{4a},\,\,t=\frac{1-q^2}{4a},\,\,u=\frac{(q^2-1)(q+1)}{4ap}.$$
\end{proposition}
This leads to the following corollary.
\begin{corollary}\label{co1}
  The most general pair $(g_K,\phi)$ on $M_l$ compatible with the $\spg(2,\bbR)$ invariant metric 
  $$g_K=(e^4)^2-2e^3\odot e^5+e^2\odot e^6+2e^1\odot e^7,$$
  coming from the Killing form in $\spg(2,\bbR)$, is a 3-parameter family with $\phi$ given by:
  $$\begin{aligned}\phi=&\frac{ap}{1-q} e^{125}+a(e^{235}-e^{127})+pe^{145}+q(e^{147}+e^{345})+\frac{p(1-q)}{4a}e^{156}+\\&\frac{1-q^2}{4a}(e^{356}-e^{167})+\frac{a(q-1)}{p}e^{237}+\tfrac12e^{246}+\frac{q^2-1}{p}e^{347}+\frac{(q^2-1)(q+1)}{4ap}e^{367}.\end{aligned}$$
  Here $a\neq 0$, $p\neq 0$, $q\neq 1$ are free parameters, and $e^{\mu\nu\rho}=e^\mu\dz e^\nu\dz e^\rho$ as before. 
  \end{corollary}
\subsection{$G_2$ structures in general} Compatible pairs $(g,\phi)$ on 7-dimensional manifolds are interesting since they give examples of $G_2$ structures \cite{bryant}. In general, a $G_2$ structure consists of a compatible pair $(g,\phi)$ of a metric $g$ and a 3-form $\phi$ on a 7-dimensional manifold $M$. It is in addition assumed
that the 3-form $\phi$ is \emph{generic}, meaning that at every point of $M$ it lies in one of the two \emph{open} orbits of the natural action of $\glg(7,\bbR)$ on 3-forms in $\bbR^7$. The simple exceptional Lie group $G_2$ appears here as the common stabilizer in $\glg(7,\bbR)$ of both $g$ and $\phi$. 

It follows (from compatibility) that the $G_2$ structures can have metrics $g$ of only two signatures: the riemannian ones and $(3,4)$ signature ones. If the signature of $g$ is riemannian, the corresponding $G_2$ structure is related to the compact real form of the simple exceptional complex Lie group $G_2$, and in the $(3,4)$ signature case the corresponding $G_2$ structure is related to the noncompact (split) real form of the complex group $G_2$. In this sense our Corollary \ref{co1} provides a 3-parameter family of split real form $G_2$ structures on $M_l$.

$G_2$ structures can be classified according to their \emph{torsion} \cite{br0,bryant}. Making the long story short we say that every $G_2$ structure $(g,\phi)$ on $M$ defines \emph{four} forms $\tau_0$, $\tau_1$, $\tau_2$ and $\tau_3$ such that 
\be\begin{aligned}
\der\phi=&\tau_0 *\phi+3 \tau_1\dz \phi+*\tau_3\\
\der*\phi=&4\tau_1\dz *\phi+\tau_2\dz \phi,\label{brtau}
\end{aligned}
\ee
where $*$ is the Hodge dual which is defined on $p$-forms $\lambda$ by
$$*\lambda(e_{\mu_1},\dots,e_{\mu_{7-p}})\,\mathrm{vol}(g)\,=\,\lambda\dz g(e_{\mu_1})\dz\dots \dz g(e_{\mu_{7-p}}),\quad\quad X\hook g(e_\mu)=g(e_\mu,X).$$ As it is visible from equations \eqref{brtau}, which we call Bryant's \cite{br0,bryant} equations in the following, each $\tau_i$, $i=0,1,2,3$, is an $i$-form on $M$. Since there is a natural action of the group $G_2$ in $\bbR^7$, which induces its action on forms, it is further required that the 3-form $\tau_3$ has values in the 27-dimensional \emph{irreducible} representation $\bigwedge^3_{27}$ of this group, the 2-form $\tau_2$ has values in the 14-dimensional \emph{irreducible} representation $\bigwedge^2_{14}$, and the 1-form has values in the 7-dimensional \emph{irreducible} representation $\bigwedge^1_{7}$. We add that the space of 1-forms $\bigwedge^1$ is irreducible, $\bigwedge^1=\bigwedge^1_7$,  and that the $G_2$ irreducible decompositions of the spaces of 2- and 3-forms look like 
$\bigwedge^2=\bigwedge^2_7\oplus\bigwedge^2_{14}$ and $\bigwedge^3=\bigwedge^3_1\oplus \bigwedge^3_7\oplus\bigwedge^3_{27}$. By the Hodge duality, the decomposition of $\bigwedge^4$ onto $G_2$ irreducible components is similar as this for $\bigwedge^3$. Here we use the convention that the lower index $i$ in $\bigwedge^p_i$ denotes the \emph{dimension} of the corresponding representation. We further mention that the 7-dimensional representations $\bigwedge^1_7$, $\bigwedge^2_7$ and $\bigwedge^3_7$ are all $G_2$ equivalent. Also, one can see that e.g. 
$\textstyle{\bigwedge^3_{27}}=\{\alpha\in\bigwedge^3\,\,\,\mathrm{s.t.}\,\,\,\alpha\dz\phi=0\,\,\&\,\,\alpha\dz*\phi=0\}$.
\subsection{All $\spg(2,\bbR)$ symmetric $G_2$ structures on $M_l$ with the metric coming from the Killing form} The below Theorem characterizes the $G_2$ structures corresponding to compatible pairs $(g_K,\phi)$ from Corollary \ref{co1}.
\begin{theorem}\label{tl}
  Let $g_K$ be the $(3,4)$ signature metric on $M_l=\spg(2,\bbR)/\slg(2,\bbR)_l$ arising as the restriction of the Killing form $K$ from $\spg(2,\bbR)$ to $M_l$,
  $$g_K=(e^4)^2-2e^3\odot e^5+e^2\odot e^6+2e^1\odot e^7.$$
  Then the most general $G_2$ structure associated with such $g_K$ is a 3-parameter family $(g_K,\phi)$ with the 3-form 
  $$\begin{aligned}\phi=&\frac{ap}{1-q} e^{125}+a(e^{235}-e^{127})+pe^{145}+q(e^{147}+e^{345})+\frac{p(1-q)}{4a}e^{156}+\\&\frac{1-q^2}{4a}(e^{356}-e^{167})+\frac{a(q-1)}{p}e^{237}+\tfrac12e^{246}+\frac{q^2-1}{p}e^{347}+\frac{(q^2-1)(q+1)}{4ap}e^{367}.\end{aligned}$$
  For this structure the torsions $\tau_\mu$ solving the Bryant's equations \eqref{brtau} are:
$$\begin{aligned}
\tau_0=&\,\frac{6}{7}\,\frac{(2a-p)^2q-(2a+p)^2}{ap},\\
\tau_1=&\frac{1}{4}\,(2a-p)\,\Big(\,-\,e^2\,+\,\frac{1}{2}\,\frac{(2a+p)(q-1)}{ap}\,e^4\,+\,\frac{1}{2}\,\frac{q^2-1}{ap}\,e^6\,\Big),\\
&\\
\tau_2=&\,0
,
\end{aligned}$$
$$\begin{aligned}
    \tau_3=&\,\Big(\tfrac{3}{28}(2a-p)^2+\frac{8ap}{7(q-1)}\Big)e^{125}+\frac{11p^2+16ap-12a^2+3q(2a-p)^2}{28p}e^{127}-\\
    &\frac{44a^2+16ap-3p^2+3q(2a-p)^2}{28a}e^{145}+\frac{(7-4q)(2a+p)^2-3q^2(2a-p)^2}{28ap}e^{147}+\\
    &\frac{3p^2(q-1)^2-12ap(q^2-1)+4a^2(31+22q+3q^2)}{112a^2}e^{156}-\frac{(q^2-1)(44a^2+16ap-3p^2+3q(2a-p)^2)}{112a^2p}e^{167}+\\
    &\frac{12a^2-16ap-11p^2-3q(2a-p)^2}{28p}e^{235}-\frac{12a^2(q-1)^2-12ap(q^2-1)+p^2(31+22q+3q^2)}{28p^2}e^{237}+\\
    &\frac{4ap(6-q)+(4a^2+p^2)(q-1)}{14ap}e^{246}+\frac{(7-4q)(2a+p)^2-3q^2(2a-p)^2}{28ap}e^{345}+\\
    &\frac{(q^2-1)(12a^2-16ap-11p^2-3q(2a-p)^2)}{28ap^2}e^{347}+\frac{(q^2-1)(44a^2+16ap-3p^2+3q(2a-p)^2)}{112a^2p}e^{356}+\\
    &\frac{(q^2-1)(q+1)(12a^2-44ap+3p^2-3q(2a-p)^2)}{112a^2p^2}e^{367},
  \end{aligned}$$
where, as usual $e^{\mu\nu}=e^\mu\dz e^\nu$ and $e^{\mu\nu\rho}=e^\mu\dz e^\nu\dz e^\rho$.  
\end{theorem}
Thus the 3-parameter family of $G_2$ structures on $M_l$ described in this Theorem have the entire 14-dimensional torsion $\tau_2=0$. This means that \emph{all} these $G_2$ structures \emph{are integrable} in the terminology of \cite{F,FI}, or what is the same, this means that they all have the \emph{totally skew symmetric torsion}.
\section{$G_2$ structures on $\spg(2,\bbR)/\slg(2,\bbR)_s$}
Now we consider the homogeneous space $M_s=\spg(2,\bbR)/\slg(2,\bbR)_s$. Since $\slg(2,\bbR)$ is spanned by $E_2,E_5+E_6,E_9$ it is convenient to put these vectors at the end of the new basis of the Lie algebra $\spg(2,\bbR)$. We choose this new basis $(f_I)$ in $\spg(2,\bbR)$ as:
$$f_1=E_1,\,\,f_2=E_3,\,\,f_3=E_4,\,\,f_4=E_6-E_5,\,\,f_5=E_7,\,\,f_6=E_8,\,\,f_7=E_{10},\,\,f_8=E_2,\,\,f_9=E_5+E_6,\,\,f_{10}=E_{9}.$$
If now, one considers $(f_I)$ as the basis of the Lie algebra of invariant vector fields on the Lie group $\spg(2,\bbR)$ then the dual basis $(f^I)$, $f_I\hook f^J=\delta^J{}_I$, of the left invariant forms on $\spg(2,\bbR)$ satisfies:
\be\begin{aligned}
  \der f^1=&2f^1\dz(f^4-f^9)+f^3\dz f^8\\
    \der f^2=&-2f^2\dz (f^4+f^9)+2f^5\dz f^8\\
      \der f^3=&2f^1\dz f^{10}+2f^3\dz f^4+f^6\dz f^8\\
        \der f^4=&2f^1\dz f^7+\tfrac12f^2\dz f^6+f^3\dz f^5\\
          \der f^5=&f^2\dz f^{10}+2f^4\dz f^5-2f^7\dz f^8\\
            \der f^6=&2f^3\dz f^{10}-2(f^4+f^9)\dz f^6\\
              \der f^7=&2(f^4-f^9)\dz f^{7}-f^5\dz f^{10}\\
                \der f^8=&2f^1\dz f^5+f^2\dz f^3-2f^8\dz f^9\\
                  \der f^9=&-2f^1\dz f^7+\tfrac12f^2\dz f^6+f^8\dz f^{10}\\
                    \der f^{10}=&2f^3\dz f^7-f^5\dz f^6-2f^9\dz f^{10}.
\end{aligned}\label{mc4}\ee
In this basis the Killing form on $\spg(2,\bbR)$ is
$$K=\tfrac{1}{12}c^I{}_{JK}c^K{}_{LI}f^J\odot f^L=2(f^4)^2-2f^3\odot f^5+f^2\odot f^6-4f^1\odot f^7+2(f^9)^2+2f^8\odot f^{10},$$
where as usual the structure constants $c^I{}_{JK}$ are defined by $[f_I,f_J]=c^K{}_{IJ}f_K$. 

Using the same arguments, as in the case of $M_l$, we again see that $\spg(2,\bbR)$ has the structure of the principal $\slg(2,\bbR)$ fiber bundle $\slg(2,\bbR)_s\to \spg(2,\bbR)\to M_s=\spg(2,\bbR)/\slg(2,\bbR)_s$ over the homogeneous space $M_s=\spg(2,\bbR)/\slg(2,\bbR)_s$. In particular we have a foliation of $\spg(2,\bbR)$ by integral leaves of an integrable distribution $D_s$ spanned by the annihilator of the forms $(f^1,f^2,\dots,f^7)$. As before, also in this Section, we will use Greek indices $\mu,\nu$, etc., to run from 1 to 7. They now number the first seven basis elements in the bases $(f_I)$ and $(f^I)$.

Repeating the procedure from the previous Sections, I now search for all bilinear symmetric forms $g=g_{\mu\nu}f^\mu\odot f^\nu$ on $\spg(2,\bbR)$, with constant coefficients $g_{\mu\nu}=g_{\nu\mu}$, whose Lie derivative with respect to any vector field $X$ from $D_l$ vanishes, \be{\mathcal L}_Xg=0\,\, \mathrm{for\,\, all}\,\, X\,\, \mathrm{in}\,\, D_s.\label{c14}\ee

I have the following proposition.  
\begin{proposition}\label{pr14}
  The most general $g=g_{\mu\nu}f^\mu\odot f^\nu$ satisfying condition \eqref{c14} is
  $$g=g_{33}\big((f^3)^2-2f^1\odot f^6\big)+g_{44}(f^4)^2+g_{55}\big((f^5)^2+2f^2\odot f^7\big)+2g_{26}\big(-2f^3\odot f^5+f^2\odot f^6-4 f^1\odot f^7\big).$$
\end{proposition}
Thus, this time, I only have a 4-parameter family of bilinear forms on $\spg(2,\bbR)$ that \emph{descend} to well defined pseudoriemannian metrics on the leaf space $M_s$. Note that, the restriction of the Killing form $K$ to the space where $(f^8,f^9,f^{10})\equiv 0$ is in this family. This corresponds to $g_{33}=g_{55}=0$ and $g_{44}=2$, $g_{26}=1/2$.

Again for simplicity reasons, I will solve the problem of finding $\spg(2,\bbR)$ invariant $G_2$ structures on $M_s$ restricting to only those pairs $(g,\phi)$ for which $g=g_K$, where
\be g_K=2(f^4)^2-2f^3\odot f^5+f^2\odot f^6-4f^1\odot f^7,\label{metl4}\ee
i.e. I only will consider one metric, the one coming from the restriction of the Killing form of $\spg(2,\bbR)$ to $M_s$. It is a well defined $(3,4)$ signature metric on the quotient space $M_s=\spg(2,\bbR)/\slg(2,\bbR)_s$.

I now look for the 3-forms $\phi=\tfrac16\phi_{\mu\nu\rho}f^\mu\dz f^\nu\dz f^\rho$ on $\spg(2,\bbR)$ which are such that
\be{\mathcal L}_X\phi=0\,\, \mathrm{for\,\, all}\,\, X\,\, \mathrm{in}\,\, D_s.\label{c24}\ee

I have the following proposition.
\begin{proposition}\label{pr24}
  There is precisely a 5-parameter family of 3-forms $\phi=\tfrac16\phi_{\mu\nu\rho}f^\mu\dz f^\nu\dz f^\rho$ on $\spg(2,\bbR)$ which satisfies condition \eqref{c24}.
  The general formula for $\phi$ is:
  $$\phi=a(4f^{147}+f^{246}+2f^{345})+b(2f^{156}+f^{236}-4f^{137})+qf^{136}+h(f^{256}-4f^{157}-2f^{237})+pf^{257}.$$
  Here $f^{\mu\nu\rho}=f^\mu\dz f^\nu\dz f^\rho$, and $a$, $b$, $q$, $h$ and $p$ are real constants.
    \end{proposition}
Solving for all 3-forms $\phi$ from this 5-parameter family that are compatible, as in \eqref{coml}, with the metric $g_K$ from \eqref{metl4}, I arrive at the following proposition.
\begin{proposition}
  The general solution to the equations \eqref{coml} is given by
  $$a=\tfrac12,\,\,b=h=0,\,\,p=\frac{1}{q}.$$
\end{proposition}
This leads to the following corollary.
\begin{corollary}\label{co14}
  The most general pair $(g_K,\phi)$ on $M_s$ compatible with the $\spg(2,\bbR)$ invariant metric 
  $$g_K=2(f^4)^2-2f^3\odot f^5+f^2\odot f^6-4f^1\odot f^7,$$
  coming from the Killing form in $\spg(2,\bbR)$, is a 1-parameter family with $\phi$ given by:
  $$\begin{aligned}\phi=&
2f^{147}+\tfrac12f^{246}+f^{345}+qf^{136}+\frac{1}{q}f^{257}
.\end{aligned}$$
  Here $q\neq 0$ is a free parameter, and $f^{\mu\nu\rho}=f^\mu\dz f^\nu\dz f^\rho$ as before. 
  \end{corollary}
\subsection{All $\spg(2,\bbR)$ symmetric $G_2$ structures on $M_s$ with the metric coming from the Killing form} Now I characterize the $G_2$ structures corresponding to compatible pairs $(g_K,\phi)$ from Corollary \ref{co14}.
\begin{theorem}\label{tl4}
  Let $g_K$ be the $(3,4)$ signature metric on $M_s=\spg(2,\bbR)/\slg(2,\bbR)_s$ arising as the restriction of the Killing form $K$ from $\spg(2,\bbR)$ to $M_s$,
  $$g_K=2(f^4)^2-2f^3\odot f^5+f^2\odot f^6-4f^1\odot f^7.$$
  Then the most general $G_2$ structure associated with such $g_K$ is a 1-parameter family $(g_K,\phi)$ with the 3-form 
  $$\phi=
2f^{147}+\tfrac12f^{246}+f^{345}+qf^{136}+\frac{1}{q}f^{257}.$$
For this structure
$$\der*\phi=0,$$
i.e. the torsions
$$\tau_1=\tau_2=0.$$
The rest of the torsions solving Bryant's equations \eqref{brtau} are:
$$\begin{aligned}
  \tau_0=&\,-\frac{18}{7},\\
  \tau_3=&\,\tfrac27\,\big(\,4f^{147}+f^{246}+2f^{345}\,\big)\,-\,\tfrac37\,\big(\,q\,f^{136}+\frac{1}{q}\,f^{257}\,\big).
  \end{aligned}$$
where, as usual $f^{\mu\nu\rho}=f^\mu\dz f^\nu\dz f^\rho$; $q\neq 0$.  
\end{theorem}
So on $M_s=\spg(2,\bbR)/\slg(2,\bbR)_s$ there exists a 1-parameter family of the above $G_2$ structures which is \emph{coclosed}. Therefore, in particular, it is \emph{integrable}

I note that formally I can also obtain coclosed $G_2$ structures on $M_l$, using the Theorem \ref{tl}. It is enough to take $p=2a$ in the solutions of this Theorem. The question if in the resulting 2-parameter family of the coclosed $G_2$ structures there is a 1-parameter subfamily equivalent to the structures I have on $M_s$ via Theorem \ref{tl4} needs further investigation. However, I doubt that the answer to this question is positive, since it is visible from the root diagram for $\spg(2,\bbR)$ that the spaces $M_l$ and $M_s$ are geometrically quite different. Indeed, apart from the $\spg(2,\bbR)$ invariant $G_2$ structures, which I have just introduced in this note, the spaces $M_l$ and $M_s$  have quite different \emph{additional} $\spg(2,\bbR)$ invariant structures.
A short look at the root diagram on page 1 of this note, shows that $M_l$ has \emph{two} well defined $\spg(2,\bbR)$ invariant rank 3-distributions, corresponding to the pushforwards from $\spg(2,\bbR)$ to $M_l$ of the vector spaces $D_{l1}=\Span_\bbR(E_2,E_3,E_7)$ and $D_{l2}=\Span_\bbR(E_4,E_8,E_9)$. Likewise $M_s$, in addition to the discussed $G_2$ structures, has also a well defined pair of $\spg(2,\bbR)$ invariant rank 3-distributions, corresponding to the pushforwards from $\spg(2,\bbR)$ to $M_s$ of the vector spaces $D_{s1}=\Span_\bbR(E_1,E_4,E_8)$ and $D_{s2}=\Span_\bbR(E_3,E_7,E_{10})$. The problem is that these two sets of pairs of $\spg(2,\bbR)$ invariant distributions are quite different. The distributions on $M_l$ have constant growth vector $(2,3)$, while the distributions on $M_s$ are \emph{integrable}. These pairs of distributions constitute an immanent ingredient of the geometry on the corresponding spaces $M_l$ and $M_s$ and, since they are diffeomorphically nonequivalent, they make the $G_2$ geometries there quite diferent. I believe that this fact makes the $G_2$ structures obtained on $M_l$ and $M_s$ really nonequivalent.    
 
\end{document}